\documentclass{article}

\usepackage{PRIMEarxiv}

\usepackage[table,xcdraw]{xcolor} 
\usepackage[utf8]{inputenc} 
\usepackage[T1]{fontenc}    
\usepackage{hyperref}       
\usepackage{url}            
\usepackage{booktabs}       
\usepackage{amsfonts}       
\usepackage{nicefrac}       
\usepackage{microtype}      
\usepackage{lipsum}
\usepackage{fancyhdr}       
\usepackage{graphicx}       
\graphicspath{{media/}}     
\usepackage{orcidlink}
\usepackage{multirow} 
\usepackage{array}    

\usepackage{amsmath}

\usepackage{pythonhighlight}

\title{EOSpython Version 0.0.11: a framework for scenario generation and a solution system for the agile earth observation satellites scheduling problem
}

\author{
  \textbf{Alex Elkjær Vasegaard} \orcidlink{0000-0001-6306-4747}\\
  Department of Materials and production \\
  Aalborg University \\
  Denmark\\
  aev@mp.aau.dk \\
  \and
  \textbf{Andreas Kühne Larsen} \orcidlink{0009-0004-3679-2254}\\
  Department of Materials and production \\
  Aalborg University \\
  Denmark\\
  ankula@mp.aau.dk \\
}

\begin{document}
\maketitle

\begin{abstract}
EOSpython is a PyPI published Python package that encompases everything within a centralized earth observation satellite scheduling system in terms of customer database setup, scenario generation, pre-processing, problem setup, scheduling solution approach, decision maker preference integration, and visualization. The package is tailored to easily configure internal parameters and contribute with other solution approaches.
\end{abstract}


\section{Introduction}
The capabilities and number of Earth Observing Satellites (EOS) have significantly increased in the last years with launch prices now well below 5000 dollars per kg to Low Earth Orbit - with the expectation of a continuous drop in prices \cite{jones2018recent, adilov2022analysis}. Additionally, EOS applications span a wide range of fields including climate monitoring, agricultural operations, urban planning, maritime surveillance, military reconnaissance, stock market analysis, emergency response planning, and an expanding array of other domains, the demand for Earth Observing Satellites is expected to diversify and continue its growth \cite{wang2020agile}. 
At the same time as we see needs for higher imaging resolution, we see a need for higher temporal resolution, meaning \textit{"the ultimate EO satellite is more EO satellites"} as a motivator for larger constellations in the future. 

Scheduling these EOSs poses a complex challenge that involves dealing with uncertainties related to factors like cloud coverage, satellite positions, and user preferences. This scheduling task presents itself as an NP-hard multi-objective optimization problem \cite{vasegaard2020multi}. Moreover, there is a requirement for near real-time execution, which has motivated researchers to develop new solution approaches to achieve satisfactory results \cite{vasquez2001logic, wei2021deep, vasegaard2020multi}.
However, due to the growing demand for satellite services, their expanding capabilities, and numbers, the need for effective scheduling solutions remains a priority. Researchers continue to work on developing improved methods to address these challenges and meet the evolving requirements on the use of EOSs - even diving into quantum compute methods \cite{stollenwerk2021agile, makarov2024quantum}.

\begin{figure}
    \centering
    \includegraphics[width=0.9\textwidth]{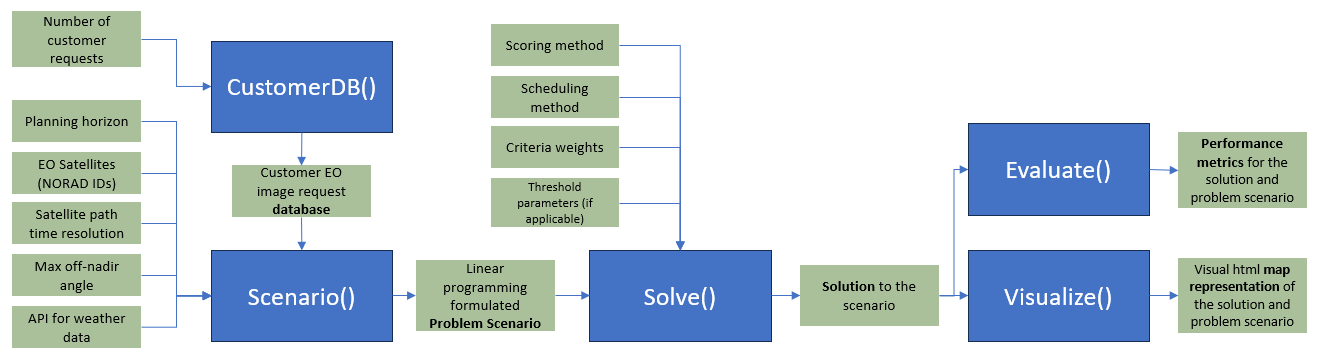}
    \caption{The functional overview of the EOSpython package framework.}
    \label{fig:functional_map}
\end{figure}

This paper introduces a PyPI published python package that encompasses everything within a centralized Earth Observation Satellite scheduling system in terms of scenario generation with customer database generation, EOS constellation selection, pre-processing, and the general optimization problem setup, then there is also a module for solution approach design and execution, encompassing even decision maker preference integration, visualization, and performance evaluation. Perhaps of highest interest to the research community is the scenario generator, that can allow research communities to build comparable benchmark problem scenarios of every size. The code is available via PyPI (\url{https://pypi.org/project/EOSpython}) and Github (\url{https://github.com/AlexVasegaard/EOS}), which also features examples of the five main functions: customerDB(), Scenario(), Solve(), Visualize(), and Evaluate(). Note the interactions can be seen in the Figure \ref{fig:functional_map}. 
The functions are based on the work of the previously published research \cite{vasegaard2020multi, vasegaard2022determining, elkjaer2021improved, vasegaard2024three}. 

\begin{figure}
    \centering
    \includegraphics[width=0.9\textwidth]{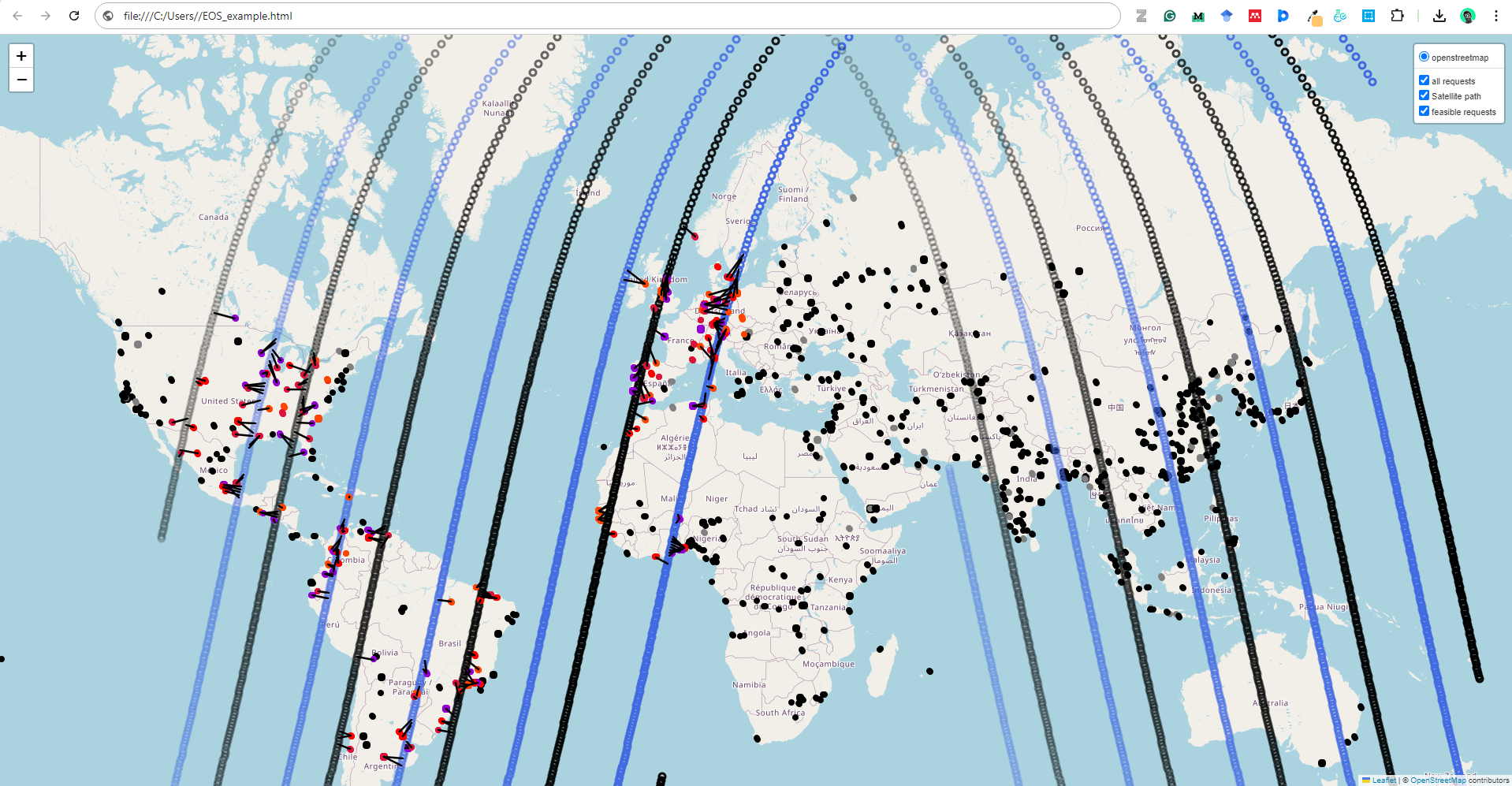}
    \caption{The interactive html output of the Visualize() function, displaying a solved scenario of the SIASP \cite{vasegaard2022determining}. The scenario is the example from Section 5.}
    \label{fig:map}
\end{figure}

\section{Scenario Generation}
The scenario generation component of the EOSpython package is (as the name implies) designed to generate EOS scheduling problems. It achieves this by generating a customer database, incorporating satellite-specific information, and creating a detailed operational scenario that includes the satellite's position, capabilities, and scheduling constraints. This section details the workflow and methods employed in generating and processing these scenarios.

\subsection{Customer Database Generation}
The first step in the scenario generation involves creating a database of customer requests, which represents the demand for satellite imaging services. The customer\_db function is responsible for this task. It synthesizes customer requests based on several parameters, including the number of requests and the satellite's swath width. This database essentially provides a list of imaging tasks that need to be assigned to available satellites.
The customer requests are distributed across the world based on population density in order to mimic real world scenarios with differing levels of saturation to the optimization problem. 

\subsection{Satellite Information and TLE Data}
With the customer database in place, the next step is to gather satellite-specific information, including the Two-Line Element (TLE) data for each satellite. The TLE data provides the necessary orbital parameters to accurately model the satellite's trajectory over time. The scenario function retrieves the TLE data from an online source using the NORAD Catalog IDs of the satellites.

Each satellite's position is computed at regular intervals within a specified time horizon, taking into account factors such as the satellite's orbital height, camera resolution, and swath width. These calculations form the basis for determining which customer requests can be serviced by each satellite at any given time and how many discrete image attempts we can compute for each request. Essentially, we discretize the continuous problem by considering some time resolution.

Other important parameters are the satellites operational imaging capabilities, that is: satellite agility (assumed linear by degrees per second), satellite swath (in km), max off-nadir angle (e.g. 30 degrees), memory capacity (in Mb), camera resolution (square meters per pixel).
Note, we assume for the memory allocation that no compressing is done on-board, and that the image stored for each acquisition is the same size that is of swath squared size. 

\subsection{Scenario Realization and Map Generation}
The last important part of the scenario realization is to define the planning horizon of the actual problem. In reality, we would like to schedule from some defined starting position to infinity, but as the satellite operations knowledge of the future becomes more and more uncertain, as well as due to practicalities of the scheduling and ability to re-schedule at a later point when communication again is feasible with the satellite, some acceptable upper bound is defined to schedule up to.  

Once the planning horizon is defined, we can deduce the satellites positions throughout. For this we have relied heavily on the python package Ephem. Here, the scenario function evaluates which customer requests can be fulfilled by calculating the distance between the satellite's ground track and the request location. This evaluation considers the satellite's operational constraints, such as its maximum off-nadir angle. We furthermore also assess the imaging characteristics, e.g. whether the sun elevation is too low for practical purposes of the imaging, or whether the forecasted cloud cover is deemed to high for consideration.

The results of this analysis, including the feasible imaging opportunities, are then visualized on an interactive map. The Folium library is used to display the satellite's path, potential acquisition points, and customer request locations. This visual output provides a clear and intuitive representation of the scenario, helping to identify the most promising imaging opportunities and potential scheduling conflicts.

\subsection{Linear Programming Problem (LPP) Setup}
The final stage of scenario generation involves setting up a Linear Programming Problem (LPP) of the satellite scheduling. The scenario function constructs a performance data frame that records all feasible imaging opportunities, considering operational imaging characteristics such as no double acquisitions, allowing for large acquisitions (multiple acquisitions of the same location), stereo images, etc. while still considering the maneuverability constraints etc. 

All this information is then used to formulate the LPP, which aims to maximize the utility of acquired images. Note, the objective function in of it self should model utility by the preferences of the satellite operation. That is, they could seek to maximize number of images acquired, number of priority 1 customer requests fulfilled, maximize the quality of the images (cloud coverage, depointing angle, sun elevation), or they could model something that considers all the different perspectives \cite{vasegaard2024three}.  

\section{Problem Formulation}
The problem has seen a variety of different formulations being presented, but by most researchers it is recognized as a generalised maximum clique problem or a generalised version of the knapsack problem \cite{vasquez2001logic, vasegaard2020multi}.

The problem formulation is based on the following assumptions related to the capabilities of the satellites:
\begin{enumerate}
    \item The preference of image attempts over other image attempts are based on a plethora of criteria and objectives, e.g. cloud coverage, profit, off-nadir angle, sun angle, customer priority, and whether the image attempt completes a stereo request or not. 
    
    \item Each satellite is equipped with a single payload capable of performing one task at a time \cite{GA_DPA_CPA_LSA}. \label{c1}
    
    \item Sequential attempts cannot always be executed due to insufficient time for satellite attitude maneuvering. The feasibility of sequential attempts depends on factors such as maneuvering time and acquisition duration. Ultimately, we assume a constant rotation speed for the satellite platform $sat_{speed}$. \label{c2}
    
    \item Requests that are too large to be acquired in one instance are divided into multiple segments, collectively representing the complete request \cite{GA_DPA_CPA_LSA}. This division is assumed to be independent from the direction of the satellite path and based on a simple heuristic rule.\label{c4}
    
    \item Within the specified time horizon, a request can be acquired only once. \label{c3}
    
    \item Two request types are considered: mono (acquired once) and stereo (acquired twice with two attempts at a convergence angle of $15^o$ to $20^o$) \cite{LS}. \label{c5}
    
    \item Satellites are subject to resource limitations in terms of memory and energy, which impact the scheduling process \cite{vasegaard2022determining}. \label{c6}
    
    \item Attempts are initiated only under specific quality thresholds, including a maximum off-nadir angle of $30^o$, a minimum sun elevation of $15^o$, and a maximum forecasted cloud cover of 60\% \cite{spot67tecnica}. \label{c7}
\end{enumerate}

On that note, the problem can be formulated as a multi-objective binary programming problem. Here one can either implement decision maker preferences through an A posterior approach as in Eq. \eqref{o1}, that is considering all the objectives in the solution search, or by A priori integration of preferences as implemented through linear scalarization technique in Eq. \eqref{o2}. Note this assumes preferences regarding objectives are additive. The objective function can consequently be formulated:
\begin{align}
    \max_x F(x) &= \max_x \{f_1(x), f_2(x),\ldots,f_n(x)\} \label{o1}\\
                &:= \max_x \sum_{i  \in N} w_{i}  x_{i} \label{o2}
\end{align}

The constraints of the problem formulation are as follows:
\begin{enumerate}
    \item Incorporating the maneuverability of the satellite in the schedule involves considering the time consumed by satellites to maneuver between pairs of attempts and the acquisition duration of the previous attempt. If the time required for changing the platforms attitude exceeds the available time (starting time of the second attempt), the maneuver is deemed infeasible.
    \begin{align}\label{maneuvarability}
        \sum_{i\in N} G_{\rho,i} x_i \leq 1 \quad \forall \rho \in \varrho
    \end{align}
    Here, the binary matrix $G$ represents infeasible pairs of attempts, serving as folded out version of the adjacency matrix for the Directed Acyclic Graph (DAG) of the problem. Here, $\varrho$ is the set of pairs of infeasible maneuvers. For each row in $G$ (representing an instance of $\rho$, we consequently, only see two 1-entries while the rest are 0, i.e. $\sum_{i \in N} G[\rho,i]=2 \phantom{0} \forall \rho \in \varrho$. Note, we can simplify $G$, by having complete sets of infeasible maneuvers. 
    \begin{align}
        G[\rho,i] = \left\{
            \begin{array}{ll}
                1 & \text{if } T_{i \rightarrow j}^{man} + T_{i}^{acq} >= t_{j}^{clock} - t_{i}^{clock} \wedge \{i,j\} \text{ is the $\rho$'th infeasible pair in $\varrho$} \\
                & \quad \text{and } \text{sat}(i) = \text{sat}(j) \\
                1 & \text{if } i=j  \wedge \{i,j\} \text{ is the $\rho$'th infeasible pair in $\varrho$}\\
                & \quad \text{and } \text{sat}(i) = \text{sat}(j) \\
                0 & \text{if } t_{j}^{clock} < t_{i}^{clock} \\
                0 & \text{otherwise}
            \end{array}
        \right.
    \end{align}
    Where $T_{i \rightarrow j}^{man}$ represents the maneuvering duration computed based on satellite platform rotational speed $sat_{speed}$ and the angle between the line of sight vectors for acquisitions $i$ and $j$ in $P$. $T_{i}^{acq}$ is the acquisition duration for attempt $i$, and $t_{i}^{clock}$ is the start time of acquisition $i$. The function $\text{sat}(.)$ returns the satellite ID of attempt $i$. This is because the maneuverability constraint only is relevant for pairs of attempts acquired by the same satellite. Additionally, we are only interested in formulating forward infeasible maneuvers, as we otherwise would have the pairs of infeasible maneuvers twice. 

    \item A request cannot be acquired multiple times unless additional acquisitions are due to the request being multi-strip or stereo. A specific number of interdependent nodes can be included in the longest path.
    \begin{align}\label{numberofrequestconstraint}
        \sum_{i\in N} B_{r,i} x_i \leq b_r \quad \forall r \in \{1,\ldots,R\}
    \end{align}
    Where $B_{r,i}$ is a binary representation of unique attempt sets $r$, and $b_r$ is the upper limit for the number of acquisitions for that unique request. Consequently, $b_r$ allows for larger image acquisitions to be divided into strips and still consider the time it takes to acquire all the strips. Note it does not state that all strips have to be acquired in the same planning horizon.
    \begin{align}
        B[r,i] = \left\{
            \begin{array}{ll}
                1 & \text{if request id $i$ is the $r$th unique id in $P$}\\
                0 & \text{otherwise}
            \end{array}
        \right.
    \end{align}

    \item Certain attempts represent stereo requests and can only be acquired together in specific pairs to fulfill a request. Only specific pairs of interdependent nodes are allowed together.
    \begin{align}\label{stereoconstraint}
        \sum_{i\in N} A[s,i] x_i = 0 \quad \forall s \in \{1,\ldots,S\}
    \end{align}
    The binary matrix $A$ represents allowed pairs $s$ of attempts for a complete stereo request. $S$ is the number of unique pairs completing a stereo request.
    \begin{align}
        A[s,i] = \left\{
            \begin{array}{ll}
                1 & \text{if attempt i is the first attempt for stereopair s} \\
                -1 & \text{if attempt i is the second attempt for stereopair s} \\
                0 & \text{otherwise}
            \end{array}
        \right.
    \end{align}
    Due to the interdependence constraint in Equation \eqref{numberofrequestconstraint}, an attempt cannot be part of more than one pair of allowed stereo requests. This also means, we have to generate a copy of the first image attempt in the stereo pair, if that image attempt is part of more stereo pairs. We then of course also have to generate a new column for every of the other constraints where this copy attempt is included. This allows us to still consider the maneuverability, while not adding too much complexity to the problem. Identifying the complexity increase is a study in of it self.

    \item Request acquisition consumes satellite memory for maneuvering. Down-link and sun-exposure considerations are neglected. Storage and energy levels are assumed to be pre-known for the scenario.
    \begin{align} 
        \sum_{i \in N} M_{i,t} x_i &\leq m_t \quad \forall t\in \Omega 
    \end{align}
    Where $M_{i,t}$ represents memory usage of attempt $i$ on satellite $t$. $\Omega$ is the set of satellites, and $m_t$ is the upper memory threshold. 
\end{enumerate}

\section{Solution Approach}
The Decision Maker preferences can be introduced through the following plethora of scoring approaches available:
\begin{enumerate}
    \item Modified ELECTRE-III 
    \item TOPSIS 
    \item Weighted Scoring Approach 
\end{enumerate}
Note, all scoring approaches contain representative variables to elicit preference information as opposed to using pairwise evaluations exhaustively. This is necessary due to the scale of the solution space. 

The solution procedure can be implemented through one of the following frameworks:
\begin{enumerate}
    \item GLPK (Large scale LPP solver from cvxopt - does not require API)
    \item ELPA (an extended longest path algorithm that can manage extremely large problem scenarios - does not require API)
    \item gurobi (large scale LPP solver from gurobi - REQUIRES API)
    \item PuLP (large scale LPP solver from PuLP - REQUIRES API)
    \item Random greedy approach (can manage large scale problems, but performs very poorly - only really applicable to showcase complexity)
\end{enumerate}
Note, some problem scenarios are too large for the LPP solvers to manage.

\section{Example on Usage}
The considered example is for the two Airbus owned satellites spot 6 and 7. Where the planning horizon is from 2024 11th of august at 10:00 to 18:00.
Where we want to acquire the set of images that maximizes the utility as much as possible. Correspondingly, the following python script can setup, model and solve the problem.

\begin{python}
from EOSpython import EOS
import pandas as pd         #The evaluation output is presented as a pd dataframe
import numpy as np          #The preference structure is given in NumPy arrays

sat_TLEs = [38755, 40053]  

horizon_start = [2024,10,17,9,40] #year,month,date,hour, minute
 
horizon = 8 #in hours

               #w
criteria_w =  [0.05,      #area
               0.1,       #off-nadir angle
               0.1,       #sun elevation
               0.2,       #cloud coverage 
               0.2,       #priority
               0.1,       #price
               0.2,       #age
               0.05]      #uncertainty

       #q,  p,   v
qpv = [[0,  30,  1000],        #area
       [0,  2,   40],          #off-nadir angle
       [0,  10,  40],          #sun elevation
       [0,  2,   15],          #cloud coverage 
       [0,  1,   4],           #priority
       [0,  100, 20000],       #price
       [0,  4,   10],          #age
       [0,  0.5,   1]]         #uncertainty

### Example with Extended Longest Path Algoorithm and ELECTRE-III scoring approach

#create database
database, map_file = EOS.customer_db(number_of_requests_0 = 300)
#Note, if map_generation is True, the database can be inspected via the interactive all_requests.html file!
 
#create scenario
x_data = EOS.scenario(customer_database = database, m = map_file,
                      seconds_gran=10, NORAD_ids=sat_TLEs, weather_real = False, 
                      schedule_start = horizon_start, hours_horizon = horizon,
                      simplify = True) #we can simplify the set of constraints when using the ELPA algorithm

x_res1 = EOS.solve(x_data, scoring_method=2, solution_method = "DAG",
                   criteria_weights_l = criteria_w, threshold_parameters_l= qpv)   #2=ELECTRE-III

EOS.visualize(x_data, x_res1, 'EOS_example') #output is an interactive map called EOS_example.html saved in the wd

df1 = EOS.evaluate(x_data, x_res1)
print(df1.solution)
print(df1.scenario)

### Example with GLPK solver and the naive weighted sum scoring approach

x_data = EOS.scenario(customer_database = database, m = map_file,
                      seconds_gran=10, NORAD_ids=sat_TLEs, weather_real = False, 
                      schedule_start = horizon_start, hours_horizon = horizon,
                      simplify = False) #for a linear solver the simplify argument must be false, as solution space otherwise is overconstrained
                      
x_res2 = EOS.solve(x_data, scoring_method=3, solution_method = "GLPK",  
                   criteria_weights_l = criteria_w, 
                   threshold_parameters_l= qpv) #3=WSA

EOS.visualize(x_data, x_res2, 'EOS_example') 
#output is an interactive map called EOS_example.html saved in the wd

df2 = EOS.evaluate(x_data, x_res2)
print(df2.solution)
print(df2.scenario)

\end{python}

\begin{table}[h!]
\centering
\rowcolors{3}{gray!25}{white} 

\begin{tabular}{|p{0.25\textwidth}|p{0.10\textwidth}|p{0.25\textwidth}|p{0.15\textwidth}|p{0.15\textwidth}|}
\hline
\multicolumn{2}{|c|}{\textbf{Problem Scenario}} & \multicolumn{3}{c|}{\textbf{Solution Performance}} \\ \hline
\textbf{Metric} & \textbf{Value} & \textbf{Metric} & \textbf{Value (ELPA)} & \textbf{Value (GLPK)} \\ \hline
Requests (\#)                   & 48        & Acquisitions (\#)         & 66     & 68       \\ \hline
Attempts (\#)                   & 472       & Total profit (euro)       & 65120  & 67598    \\ \hline
Constraints (\#)                & 842 (1125)& Avg cloud cover (\%)      & 14.7   & 14.98    \\ \hline
Avg. angle (\%)                 & 20.5      & Cloud cover < 10\%  (\#)  & 12     & 12       \\ \hline
Avg. area (sqkm)                & 553.0     & Cloud cover > 30\% (\#)   & 0      & 0     \\ \hline
Avg. price (euro)               & 994.5     & Avg. angle (degrees)      & 19.65  & 19.16     \\ \hline
Avg. sun elevation (degrees)    & 44.6      & Angle < 10\% (\#)         & 12     & 16     \\ \hline
Avg. cloud cover (\%)           & 17.2      & Angle > 30\% (\#)         & 4      & 3     \\ \hline
Avg. priority                   & 2.44      & Avg priority              & 2.58   & 2.56     \\ \hline
                                &           & Priority 1 (\#)           & 12     & 12       \\ \hline
                                &           & Priority 2 (\#)           & 19     & 21       \\ \hline
                                &           & Priority 3 (\#)           & 20     & 20       \\ \hline
                                &           & Priority 4  (\#)          & 15     & 15       \\ \hline
                                &           & Avg sun elevation (degrees)& 44.7   & 44.2 \\ \hline 
                                &           & Total area (sqkm)          & 41296  & 43124  \\
\end{tabular}
\caption{A merged performance evaluation matrix showcasing the problem scenario specification and the performance metrics for the two solution approaches in the example. Note, the number of constraints is higher for the exact approach due to the complete sets of infeasible attempts over-constraining the solution space, whereas the ELPA iteratively goes through the constraints and consequently allows for complete sets.}
\end{table}

We see by the performance evaluation matrix that we can acquire 66 or 68 images depending on the solution approach with decent off-nadir angle, sun elevation, age distribution, and cloud coverage. As well as servicing 12 high priority customers. The output can also be inspected visually by the html file portrayed in Figure\ref{fig:map}. 
Both of the solution approaches computed under 2 seconds with 1.83 and 1.53 seconds for the ELPA and GLPK, respectively. The objective function values are not to be compared directly in this example as two different objective functions are used - one based on the ELECTRE-III scoring approach given weights and threshold values \cite{vasegaard2020multi} and another based on the regular weighted sum approach given only weights to aggregate the objectives. 

Note, the relatively small scale of the problem scenario with an initial customer database of 300 is because of the exact solver GLPK's inability to solve large scale problems. The ELPA approach can easily solve problems with an initial customer database in the many thousands\cite{vasegaard2024three}. 

\section{How to contribute and collaborate}
Contributions to the package is strongly encouraged, and can be accomplished through pull-requests on Github. After a review, the changes will either be accepted or rejected. However, as the project is still it its infancy, expect variable names, function calls, and the overall structure to mutate rapidly with development. We will roll out a significant update very soon and accompanying it - a comprehensive guide on how to contribute and collaborate.

\section{Future work}


The platform agility is currently only considering agile platforms with three degrees of freedom. Additionally, the agility speed is assumed to be linear which is a stark simplification of real-world systems. 

A wide array of solution approaches should also be implemented, such as the family of evolutionary algorithms (particle-swarm optimization, etc.) and multi-objective a posterior approaches (e.g., NSGA-II). Similarly, it is a big goal that customizable objective functions could be implemented, consequently allowing for a larger (or smaller) portion of objectives and criteria to be evaluated in the final scheduling. 

A significant contribution of the EOSpython framework is the ability to add stereo imaging capabilities to the satellites. However, the novelty and proficiency of the utilized method are not very rigidly compared to other approaches, and especially from a computational efficiency perspective, this is very important. Therefore, we also expect to make further contributions to this. 

Currently, the cloud cover data is retrieved from one of two sources: either the OpenWeatherMap API or a simple cloud cover generator that takes location and time into perspective. Both methods should be improved. From the OWM perspective, the forecasts are given by relatively infrequent point forecasts, and they do not reveal the probabilistic forecast distributions that would be of great value to the EOS scheduling case. 

Moreover, we need to consider the temporal aspects of the scheduling, that is, potentially acquiring images outside the considered planning horizon and essentially saving an image attempt for later. Currently, we specify a time horizon and optimize the short-term utility within that time frame. Ideally, we want to optimize for long-term utility. 

Furthermore, the EOSpython framework allows us to create large amounts of solution and problem scenario data, which can consequently be used for training. Therefore, this work is expected to contribute to learning-based approaches in EOS scheduling in the future. 

Lastly, we expect to soon publish specific benchmark problems on which different state-of-the-art solution approaches can be compared. 

\section{Conclusion}
This paper has introduced a easy to use, customizable, and versatile EOS scheduling Python package, EOSpython, designed to address the challenges of EOS system analysis and optimization. With an array of innovative solution methods presented herein, the package serves as a valuable tool for researchers and engineers alike, facilitating the exploration and enhancement of various satellite system specifications and operational scheduling solution designs. 

\section*{Acknowledgments}
The work would not have been what it is without the inputs and sparing with Prof. Peter Nielsen, Ph.D. Mathieu Picard, and Prof. Subrata Saha. This was supported in part by Airbus Defence \& Space. 



\end{document}